\documentclass[12pt]{article}
\usepackage{latexsym, amssymb}
\usepackage{dsfont}

\DeclareMathAlphabet\gothic{U}{euf}{m}{n}

\makeatletter
\def\eqnarray{\stepcounter{equation}\let\@currentlabel=\theequation
\global\@eqnswtrue
\tabskip\@centering\let\\=\@eqncr
$$\halign to \displaywidth\bgroup\hfil\global\@eqcnt\z@
  $\displaystyle\tabskip\z@{##}$&\global\@eqcnt\@ne
  \hfil$\displaystyle{{}##{}}$\hfil
  &\global\@eqcnt\tw@ $\displaystyle{##}$\hfil
  \tabskip\@centering&\llap{##}\tabskip\z@\cr}

\def\endeqnarray{\@@eqncr\egroup
      \global\advance\c@equation\m@ne$$\global\@ignoretrue}

\def\@yeqncr{\@ifnextchar [{\@xeqncr}{\@xeqncr[5pt]}}
\makeatother

\begin{document}
\bibliographystyle{tom}

\newtheorem{lemma}{Lemma}[section]
\newtheorem{thm}[lemma]{Theorem}
\newtheorem{cor}[lemma]{Corollary}
\newtheorem{voorb}[lemma]{Example}
\newtheorem{rem}[lemma]{Remark}
\newtheorem{prop}[lemma]{Proposition}
\newtheorem{stat}[lemma]{{\hspace{-5pt}}}

\newenvironment{remarkn}{\begin{rem} \rm}{\end{rem}}
\newenvironment{exam}{\begin{voorb} \rm}{\end{voorb}}

\newcommand{\gota}{\gothic{a}}
\newcommand{\gotb}{\gothic{b}}
\newcommand{\gotc}{\gothic{c}}
\newcommand{\gote}{\gothic{e}}
\newcommand{\gotf}{\gothic{f}}
\newcommand{\gotg}{\gothic{g}}
\newcommand{\gothh}{\gothic{h}}
\newcommand{\gotk}{\gothic{k}}
\newcommand{\gotm}{\gothic{m}}
\newcommand{\gotn}{\gothic{n}}
\newcommand{\gotp}{\gothic{p}}
\newcommand{\gotq}{\gothic{q}}
\newcommand{\gotr}{\gothic{r}}
\newcommand{\gots}{\gothic{s}}
\newcommand{\gotu}{\gothic{u}}
\newcommand{\gotv}{\gothic{v}}
\newcommand{\gotw}{\gothic{w}}
\newcommand{\gotz}{\gothic{z}}
\newcommand{\gotA}{\gothic{A}}
\newcommand{\gotB}{\gothic{B}}
\newcommand{\gotG}{\gothic{G}}
\newcommand{\gotL}{\gothic{L}}
\newcommand{\gotS}{\gothic{S}}
\newcommand{\gotT}{\gothic{T}}

\newcounter{teller}
\renewcommand{\theteller}{\Roman{teller}}
\newenvironment{tabel}{\begin{list}%
{\rm \bf \Roman{teller}.\hfill}{\usecounter{teller} \leftmargin=1.1cm
\labelwidth=1.1cm \labelsep=0cm \parsep=0cm}
                      }{\end{list}}

\newcounter{tellerr}
\renewcommand{\thetellerr}{(\roman{tellerr})}
\newenvironment{subtabel}{\begin{list}%
{\rm  (\roman{tellerr})\hfill}{\usecounter{tellerr} \leftmargin=1.1cm
\labelwidth=1.1cm \labelsep=0cm \parsep=0cm}
                         }{\end{list}}

\newcommand{\Ni}{{\bf N}}
\newcommand{\Ri}{{\bf R}}
\newcommand{\Ci}{{\bf C}}
\newcommand{\Ti}{{\bf T}}
\newcommand{\Zi}{{\bf Z}}
\newcommand{\Fi}{{\bf F}}

\newcommand{\proof}{\mbox{\bf Proof} \hspace{5pt}} 
\newcommand{\remark}{\mbox{\bf Remark} \hspace{5pt}}
\newcommand{\ruimte}{\vskip10.0pt plus 4.0pt minus 6.0pt}

\newcommand{\simh}{{\stackrel{{\rm cap}}{\sim}}}
\newcommand{\ad}{{\mathop{\rm ad}}}
\newcommand{\Ad}{{\mathop{\rm Ad}}}
\newcommand{\Aut}{\mathop{\rm Aut}}
\newcommand{\arccot}{\mathop{\rm arccot}}
\newcommand{\capp}{{\mathop{\rm cap}}}
\newcommand{\rcapp}{{\mathop{\rm rcap}}}
\newcommand{\diam}{\mathop{\rm diam}}
\newcommand{\divv}{\mathop{\rm div}}
\newcommand{\codim}{\mathop{\rm codim}}
\newcommand{\RRe}{\mathop{\rm Re}}
\newcommand{\IIm}{\mathop{\rm Im}}
\newcommand{\Tr}{{\mathop{\rm Tr}}}
\newcommand{\Vol}{{\mathop{\rm Vol}}}
\newcommand{\card}{{\mathop{\rm card}}}
\newcommand{\supp}{\mathop{\rm supp}}
\newcommand{\sgn}{\mathop{\rm sgn}}
\newcommand{\essinf}{\mathop{\rm ess\,inf}}
\newcommand{\esssup}{\mathop{\rm ess\,sup}}
\newcommand{\Int}{\mathop{\rm Int}}
\newcommand{\Leibniz}{\mathop{\rm Leibniz}}
\newcommand{\lcm}{\mathop{\rm lcm}}
\newcommand{\loc}{{\rm loc}}

\newcommand{\mod}{\mathop{\rm mod}}
\newcommand{\spann}{\mathop{\rm span}}
\newcommand{\one}{\mathds{1}}

\hyphenation{groups}
\hyphenation{unitary}

\newcommand{\tfrac}[2]{{\textstyle \frac{#1}{#2}}}

\newcommand{\cb}{{\cal B}}
\newcommand{\cc}{{\cal C}}
\newcommand{\cd}{{\cal D}}
\newcommand{\ce}{{\cal E}}
\newcommand{\cf}{{\cal F}}
\newcommand{\ch}{{\cal H}}
\newcommand{\ci}{{\cal I}}
\newcommand{\ck}{{\cal K}}
\newcommand{\cl}{{\cal L}}
\newcommand{\cm}{{\cal M}}
\newcommand{\co}{{\cal O}}
\newcommand{\cs}{{\cal S}}
\newcommand{\ct}{{\cal T}}
\newcommand{\cx}{{\cal X}}
\newcommand{\cy}{{\cal Y}}
\newcommand{\cz}{{\cal Z}}

\newlength{\hightcharacter}
\newlength{\widthcharacter}
\newcommand{\covsup}[1]{\settowidth{\widthcharacter}{$#1$}\addtolength{\widthcharacter}{-0.15em}\settoheight{\hightcharacter}{$#1$}\addtolength{\hightcharacter}{0.1ex}#1\raisebox{\hightcharacter}[0pt][0pt]{\makebox[0pt]{\hspace{-\widthcharacter}$\scriptstyle\circ$}}}
\newcommand{\cov}[1]{\settowidth{\widthcharacter}{$#1$}\addtolength{\widthcharacter}{-0.15em}\settoheight{\hightcharacter}{$#1$}\addtolength{\hightcharacter}{0.1ex}#1\raisebox{\hightcharacter}{\makebox[0pt]{\hspace{-\widthcharacter}$\scriptstyle\circ$}}}
\newcommand{\scov}[1]{\settowidth{\widthcharacter}{$#1$}\addtolength{\widthcharacter}{-0.15em}\settoheight{\hightcharacter}{$#1$}\addtolength{\hightcharacter}{0.1ex}#1\raisebox{0.7\hightcharacter}{\makebox[0pt]{\hspace{-\widthcharacter}$\scriptstyle\circ$}}}

\thispagestyle{empty}

\vspace*{1cm}
\begin{center}
{\Large\bf Diffusion determines the compact manifold} \\[5mm]

\large W. Arendt$^1$ and A.F.M. ter Elst$^2$

\end{center}

\vspace{5mm}

\begin{center}
{\bf Abstract}
\end{center}

\begin{list}{}{\leftmargin=1.8cm \rightmargin=1.8cm \listparindent=10mm 
   \parsep=0pt}
\item
We provide a short proof for the theorem that two compact
Riemannian manifolds are isomorphic if and only 
there exists an order isomorphism which intertwines between the 
heat semigroups on the manifolds.

\end{list}

\vspace{4cm}
\noindent
March 2011.

\vspace{5mm}
\noindent
AMS Subject Classification: 58J53, 35P05, 47F05, 35R30.

\vspace{15mm}

\noindent
{\bf Home institutions:}    \\[3mm]
\begin{tabular}{@{}cl@{\hspace{10mm}}cl}
1. & Abteilung Angewandte Analysis & 
2. & Department of Mathematics  \\
& Universit\"at Ulm & 
  & University of Auckland  \\
& Helmholtzstr.\ 18  &
  & Private bag 92019 \\
& 89069 Ulm  &
   & Auckland 1142 \\ 
& Germany  &
  & New Zealand \\[8mm]
\end{tabular}

\newpage
\setcounter{page}{1}

\section{Introduction} \label{Scdrum1}

Two years before the publication of Kac's famous paper \cite{Kac}
`Can one hear the shape of a drum' Milnor \cite{Milnor} gave a 
counter example showing that one cannot hear the shape of a compact 
Riemannian manifold.
Milnor presented two 16-dimensional Riemannian manifolds
for which the associated Laplace--Beltrami operators have the 
same spectrum, i.e.\ are isospectral.
The latter is equivalent with the existence of a unitary operator $U$
which intertwines the heat semigroups on the compact manifolds.
The heat semigroups are positive, which means that they map positive functions
(i.e.\ positive heat) to positive functions on the $L_2$-spaces 
of the compact manifolds.
In this paper we replace the unitary operator by an 
order isomorphism, i.e.\ a linear bijective mapping $U$ such that 
$U \varphi \geq 0$ if and only if $\varphi \geq 0$.
Then we show that the manifolds are indeed isomorphic.
This may be interpreted in the following way.
The heat semigroups are positive, which means that positive functions
(heat densities) are mapped to positive functions.
The orbit corresponding to a positive initial value describes the 
propagation of the heat density, i.e.\ the diffusion.
Thus to say that an order isomorphism intertwines between two heat 
semigroups means that the positive orbits are mapped to positive 
orbits.
So our result may be rephrased by saying that diffusion determines 
the compact manifold.
For open connected subsets of $\Ri^d$ satisfying a weak smoothness
condition Arendt \cite{Are3} proved that diffusion determines the body
(see also \cite{Are4}).
In a recent paper \cite{ABE} this was extended to connected Riemannian 
manifolds satisfying the same smoothness condition.
Every compact connected Riemannian manifold satisfies this smoothness condition.

The aim of this paper is to give a direct and short proof that 
diffusion determines the body for compact  Riemannian manifolds.
The compact Riemannian manifolds do not have to be connected.

\smallskip

Let $(M,g)$ be a compact Riemannian manifold of dimension $d$.
Then $M$ has a natural Radon measure with respect to which we define the 
$L_p$-spaces on $M$. 
Set 
\[
H^1(M)
= \{ \varphi \in L_2(M) : \varphi \circ x^{-1} \in H^1(x(V))
       \mbox{ for every chart } (V,x) \} 
\;\;\; .  \]
If $\varphi \in H^1(M)$ and $(V,x)$ is a chart on $M$ then set 
$\frac{\partial}{\partial x^i} \varphi = (D_i (\varphi \circ x^{-1})) \circ x \in L_2(V)$,
where $D_i$ denotes the partial derivative in $\Ri^d$.
Moreover, for all $\varphi,\psi \in H^1(M)$
there exists a unique element $\nabla \varphi \cdot \nabla \psi \in L_1(M)$
such that 
\[
\nabla \varphi \cdot \nabla \psi \Big|_V
= \sum_{i,j=1}^d g^{ij} \Big( \frac{\partial}{\partial x^i} \varphi \Big)
                               \Big( \frac{\partial}{\partial x^j} \psi \Big)
\]
for every chart $(V,x)$ on $M$.
Set $|\nabla \varphi| = (\nabla \varphi \cdot \nabla \varphi)^{1/2}$.
We provide $H^1(M)$ with the norm
$\varphi \mapsto ( \|\varphi\|_2^2 + \| \, |\nabla \varphi| \, \|_2^2 )^{1/2}$.
Then $H^1(M)$ is a Hilbert space.
Define the bilinear form $a \colon H^1(M) \times H^1(M) \to \Ri$ by
$a(\psi,\varphi) = \int \nabla \psi \cdot \nabla \varphi$.
Then $a$ is a closed and positive form in $L_2(M)$.
The {\bf Dirichlet Laplace--Beltrami operator} $\Delta$
on $M$ is the associated self-adjoint operator.
If $(V,x)$ is a chart on $M$ then 
\[
\Delta \, \varphi
= - \sum_{i,j=1}^d \frac{1}{\sqrt{g}} \, \frac{\partial}{\partial x^i} \, 
           g^{ij} \, \sqrt{g} \, \frac{\partial}{\partial x^j} \, \varphi
\]
for all $\varphi \in C_c^\infty(V)$.
Let $S$ be the semigroup on $L_2(M)$ generated by $-\Delta$
and let $p \in [1,\infty)$.
By the Beurling--Deny criteria the operator
$S_t|_{L_2(M) \cap L_p(M)}$ extends to a positive contraction operator on 
$L_p(M)$ for all $t > 0$.
Moreover, $S^{(p)}$ is a $C_0$-semigroup.
Since the semigroup $S$ has a smooth kernel satisfying
Gaussian bounds (\cite{Sal5} Theorem~5.4.12), it follows that 
$S_t C(M) \subset C(M)$ and
$S|_{C(M)}$ is a $C_0$-semigroup on $C(M)$.

If $(M_1,g_1)$ and $(M_2,g_2)$ are two compact Riemannian manifolds 
then a map $\tau \colon M_1 \to M_2$ is called an {\bf isometry}
if it is a $C^\infty$-diffeomorphism and 
\[
g_2|_{\tau(p)}(\tau_*(v), \tau_*(w)) = g_1|_p(v,w)
\]
for all $p \in M_1$ and $v,w \in T_p M_1$.
The Riemannian manifolds $(M_1,g_1)$ and $(M_2,g_2)$
are called {\bf isomorphic} if there exists an isometry 
from $M_1$ onto $M_2$.
If $\tau \colon M_1 \to M_2$ is an isometry
and $p \in [1,\infty]$ then 
$\varphi \circ \tau \in L_p(M_1)$ and
\begin{equation}
\|\varphi \circ \tau\|_{L_p(M_1)} = \|\varphi\|_{L_p(M_2)}
\label{eSdrum1;1}
\end{equation}
for all $\varphi \in L_p(M_2)$.

A linear operator $U \colon E \to F$ between two Riesz spaces is said to be a 
{\bf lattice homomorphism} if 
\[
U (\varphi \wedge \psi) = (U \varphi) \wedge (U \psi)
\]
for all $\varphi,\psi \in E$.
For alternative equivalent definitions see \cite{AB} Theorem~7.2.
Each lattice homomorphism $U$ is positive, i.e.\ $\varphi \geq 0$
implies $U \varphi \geq 0$.
An {\bf order isomorphism} $U \colon E \to F$ is a bijective 
mapping  such that $U \varphi \geq 0$ if and only if $\varphi \geq 0$.
Equivalently, $U$ is an order isomorphism
if and only if $U$ is a bijective lattice homomorphism.
Then also $U^{-1}$ is an order isomorphism.
Recall also that each positive operator between $L_p$-spaces,
or from $C(M_1)$ into $C(M_2)$ where $M_1$ and $M_2$ are 
compact Hausdorff spaces,  is continuous 
by \cite{AB} Theorem~12.3.

The main theorem of this paper is the following.

\begin{thm} \label{tcdrum101}
Let $(M_1,g_1)$ and $(M_2,g_2)$ be two compact Riemannian manifolds.
Let $p \in [1,\infty)$.
For all $j \in \{ 1,2 \} $ let $\Delta_j$ be the Laplace--Beltrami operator 
on $M_j$ and let $S^{(j)}$ and $T^{(j)}$ be the associated semigroups on $L_p(M_j)$
and $C(M_j)$.
Then the following three conditions are equivalent.
\begin{tabel}
\item \label{tcdrum101-1}
$(M_1,g_1)$ and $(M_2,g_2)$ are isomorphic.
\item \label{tcdrum101-3}
There exists an order isomorphism $U \colon L_p(M_1) \to L_p(M_2)$ such that 
\[
U S^{(1)}_t = S^{(2)}_t U
\]
for all $t > 0$.
\item \label{tcdrum101-4}
There exists an order isomorphism $U \colon C(M_1) \to C(M_2)$ such that 
\[
U T^{(1)}_t = T^{(2)}_t U
\]
for all $t > 0$.
\end{tabel}
Moreover, if the manifolds are connected and if
$U$ is an order isomorphism as in 
Condition~{\rm \ref{tcdrum101-3}} or {\rm \ref{tcdrum101-4}} then 
there exist $c>0$ and a $($surjective$)$ isometry $\tau \colon M_2 \to M_1$ such that 
$U \varphi = c \, \varphi \circ \tau$ for all $\varphi \in L_p(M_1)$.
\end{thm}

The implications \ref{tcdrum101-1}$\Rightarrow$\ref{tcdrum101-3} and 
\ref{tcdrum101-1}$\Rightarrow$\ref{tcdrum101-4} are an easy consequence of 
(\ref{eSdrum1;1}).

\section{Proof of Theorem~\ref{tcdrum101}}

The first part in the proof of Theorem~\ref{tcdrum101} is the 
observation that $C^\infty$-functions are invariant under intertwining 
operators.

\begin{lemma} \label{lcdrum201}
Let $(M_1,g_1)$ and $(M_2,g_2)$ be two compact Riemannian manifolds.
Let $p \in [1,\infty)$.
For all $j \in \{ 1,2 \} $ let $\Delta_j$ be the Laplace--Beltrami operator 
on $M_j$ and let $S^{(j)}$ and $T^{(j)}$ be the associated semigroups on $L_p(M_j)$
and $C(M_j)$.
Let either $U \colon L_p(M_1) \to L_p(M_2)$ be an order isomorphism such that 
\[
U S^{(1)}_t = S^{(2)}_t U
\]
for all $t > 0$, or
$U \colon C(M_1) \to C(M_2)$ be an order isomorphism such that 
\begin{equation}
U T^{(1)}_t = T^{(2)}_t U
\label{elcdrum201;2}
\end{equation}
for all $t > 0$.
Then
\begin{subtabel}
\item \label{lcdrum201-1}
$U C^\infty(M_1) = C^\infty(M_2)$.
\item \label{lcdrum201-2}
$U \varphi \geq 0$ if and only if  $\varphi \geq 0$, for all $\varphi \in C^\infty(M_1)$. 
\item \label{lcdrum201-3}
$(U \varphi) (U \psi) = 0$ for all $\varphi,\psi \in C^\infty(M_1)$ with $\varphi \, \psi = 0$.
\item \label{lcdrum201-4}
$\Delta_2 U \varphi = U \Delta_1 \varphi$ for all $\varphi \in C^\infty(M_1)$.
\end{subtabel}
\end{lemma}
\proof\
Suppose $U$ is an order isomorphism from $C(M_1)$ onto $C(M_2)$.
Let $H_j$ be the generator of $T^{(j)}$ for all $j \in \{ 1,2 \} $.
If $\varphi \in D(H_1)$ then it follows from (\ref{elcdrum201;2}) that 
\[
\tfrac{1}{t} (I - T_t^{(2)}) U \varphi
= \tfrac{1}{t} U (I - T^{(1)}) \varphi
\]
for all $t > 0$.
Since $U$ is continuous one deduces that 
$U \varphi \in D(H^{(2)})$.
So $U D(H_1) \subset D(H_2)$ and $H_2 U \varphi = U H_1 \varphi$
for all $\varphi \in D(\Delta_1)$.
Similarly $U^{-1} D(H_2) \subset D(H_1)$ and therefore $U D(H_1) = D(H_2)$.
Hence by iteration
$U \bigcap_{n=1}^\infty D(H_1^n) = \bigcap_{n=1}^\infty D(H_2^n)$.
But $C^\infty(M_j) = \bigcap_{n=1}^\infty D(H_j^n)$ 
for all $j \in \{ 1,2 \} $ by elliptic regularity.
Here we use that the manifolds are compact.
This shows \ref{lcdrum201-1} and \ref{lcdrum201-4}.
Property~\ref{lcdrum201-2} follows since $U$ is an order isomorphism.
Moreover, 
$|U \varphi| = U |\varphi|$ for all $\varphi \in C(M_1)$.
Hence if $\varphi,\psi \in C(M_1)$ and $\varphi \, \psi = 0$ then 
$|\varphi| \wedge |\psi| = 0$ and 
$|U \varphi| \wedge |U \psi| = U |\varphi| \wedge U |\psi|
= U(|\varphi| \wedge |\psi|) = 0$.
Therefore $|(U \varphi) (U \psi)| = |U \varphi| \, |U \psi| = 0$
and $(U \varphi) (U \psi) = 0$.
This implies Property~\ref{lcdrum201-3}.

The proof on the $L_p$-spaces is similar.\hfill$\Box$

\ruimte

The next lemma is a $C^\infty$-version of the Riesz representation theorem.
(Cf.\ \cite{EvG} Corollary~1.8.1.)

\begin{lemma} \label{lcdrum201.3}
Let $M$ be a compact Riemannian manifold and $F \colon C^\infty(M) \to \Ri$ a 
positive linear functional
such that 
\begin{equation}
F(\varphi) \, F(\psi) = 0
\mbox{ for all } \varphi,\psi \in C^\infty(M) \mbox{ with } \varphi \, \psi = 0
.  
\label{elcdrum201.3;1}
\end{equation}
Then there exist $c \in [0,\infty)$ and $p \in M$ such that 
$F(\varphi) = c \, \varphi(p)$ for all $\varphi \in C^\infty(M)$.
\end{lemma}
\proof\
Let $\varphi \in C^\infty(M)$.
Then $\|\varphi\|_\infty \, \one - \varphi \geq 0$, so it follows from positivity
that $F(\varphi) \leq F(\one) \, \|\varphi\|_\infty$.
Since $C^\infty(M)$ is dense in $C(M)$ one can extend $F$ 
to a continuous linear function from $C(M)$ into $\Ri$.
This extension is again positive since positive functions in $C(M)$ 
can be approximated uniformly by positive functions in $C^\infty(M)$.
By the Riesz representation theorem
there exists a unique Radon measure $\mu$ on $M$ such that
$F(\varphi) = \int \varphi \, d \mu$ for all $\varphi \in C^\infty(M)$.
Then it follows from (\ref{elcdrum201.3;1}) that $\mu$ is a point measure.
Hence there exist $p \in M$ and $c \in [0,\infty)$ such that 
$F(\varphi) = c \, \varphi(p)$ for all $\varphi \in C^\infty(M)$.\hfill$\Box$

\begin{prop} \label{pdrum202}
Let $(M_1,g_1)$ and $(M_2,g_2)$ be two compact Riemannian manifolds.
Suppose there exists a linear bijection $U \colon C^\infty(M_1) \to C^\infty(M_2)$
such that 
\begin{subtabel}
\item \label{pdrum202-1}
$U \varphi \geq 0$ if and only if  $\varphi \geq 0$, for all $\varphi \in C^\infty(M_1)$. 
\item \label{pdrum202-2}
$(U \varphi) (U \psi) = 0$ if and only if $\varphi \, \psi = 0$, for all $\varphi,\psi \in C^\infty(M_1)$.
\item \label{pdrum202-3}
$\Delta_2 U \varphi = U \Delta_1 \varphi$ for all $\varphi \in C^\infty(M_1)$.
\end{subtabel}
Then the Riemannian manifolds $(M_1,g_1)$ and $(M_2,g_2)$ are isomorphic.
\end{prop}
\proof\
Let $q \in M_2$.
Then the map $\varphi \mapsto (U \varphi)(q)$ from $C^\infty(M_1)$ into 
$\Ri$ is linear, positive and non-zero.
So by Lemma~\ref{lcdrum201.3} there exist $\tau(q) \in M_1$ and $h(q) \in (0,\infty)$ such that 
\begin{equation}
(U \varphi)(q)
= h(q) \, \varphi(\tau(q))
\label{epdrum202;1}
\end{equation}
for all $\varphi \in C^\infty(M_1)$.
So one obtains functions $\tau \colon M_2 \to M_1$ and 
$h \colon M_2 \to (0,\infty)$.
Similarly, there exist $\widetilde \tau \colon M_1 \to M_2$
and $\tilde h \colon M_1 \to (0,\infty)$ such that 
$(U^{-1} \psi)(p) = \tilde h(p) \, \psi(\widetilde \tau(p))$
for all $\psi \in C^\infty(M_2)$ and $p \in M_1$.
Then 
$\varphi(p) 
= \tilde h(p) \, h(\widetilde \tau(p)) \, \varphi(\tau( \widetilde \tau(p)))$
for all $\varphi \in C^\infty(M_1)$ and $p \in M_1$.
Choosing $\varphi = \one$ gives $\tilde h(p) \, h(\widetilde \tau(p)) = 1$.
Hence $\varphi = \varphi \circ \tau \circ \widetilde \tau$
for all $\varphi \in C^\infty(M_1)$ and
$\tau \circ \widetilde \tau = I$.
Similarly $\widetilde \tau \circ \tau = I$ and $\tau$ is a bijection.

Choosing again $\varphi = \one$ in (\ref{epdrum202;1}) gives
$h = U \one \in C^\infty(M_2)$.
Hence $\varphi \circ \tau = h^{-1} U \varphi \in C^\infty(M_2)$ for all 
$\varphi \in C^\infty(M_1)$ and $\tau$ is a $C^\infty$-function.
Thus $\tau$ is a $C^\infty$-diffeomorphism and the two manifolds have 
the same dimension.
Let $d = \dim M_1 = \dim M_2$.

It follows from Property~\ref{pdrum202-3} that 
\begin{equation}
\Delta_2( h \cdot (\varphi \circ \tau))
= h \cdot ( (\Delta_1 \varphi) \circ \tau)
\label{epdrum202;2}
\end{equation}
for all $\varphi \in C^\infty(M_1)$.
Let $q \in M_2$.
There exists a chart $(V,x)$ on $M_1$ such that 
$\tau(q) \in V$ and $x(\tau(q)) = 0$.
Let $\Omega \subset M_1$ be open such that 
$\tau(q) \in \Omega \subset \overline\Omega \subset V$.
Let $\lambda_1,\ldots,\lambda_d \in \Ri$.
For all $t > 0$ there exists a $\varphi_t  \in C^\infty(M_1)$
such that 
\[
\varphi_t |_\Omega 
= e^{t \sum_{k=1}^d \lambda_k x^k} |_\Omega
.  \]
Since
\[
\Delta_1
= \sum_{i,j=1}^d \frac{1}{\sqrt{g_1}} \, \frac{\partial}{\partial x^i} \, 
           g_1^{ij} \, \sqrt{g_1} \, \frac{\partial}{\partial x^j}
\]
on $V$ it follows that 
\[
\Delta_1 \varphi_t
= \sum_{i,j=1}^d t^2 g_1^{ij} \, \lambda_i \, \lambda_j \, \varphi_t
   - t \, \frac{\lambda_j}{\sqrt{g_1}} \, \varphi_t \, 
     \frac{\partial}{\partial x^i} ( g_1^{ij} \, \sqrt{g_1} )
\]
on $\Omega$.
Hence 
\[
\lim_{t \to \infty} t^{-2} \Big( h \cdot ( (\Delta_1 \varphi_t) \circ \tau) \Big) (q)
= h(q) \sum_{i,j=1}^d g_1^{ij}(\tau(q)) \, \lambda_i \, \lambda_j
.  \]
Next, $(\tau^{-1}(V),y)$ is a chart on $M_2$, where $y = x \circ \tau$.
Then it follows similarly that 
\begin{eqnarray*}
\lim_{t \to \infty} t^{-2} \Big( \Delta_2( h \cdot (\varphi_t \circ \tau)) \Big) (q)
& = & \sum_{i,j=1}^d h(q) \, g_2^{ij}(q)
   \Big( \frac{\partial}{\partial y_i} \sum_{k=1}^d \lambda_k x^k \circ \tau \Big)(q)
   \Big( \frac{\partial}{\partial y_j} \sum_{l=1}^d \lambda_l x^l \circ \tau \Big)(q)  \\
& = & \sum_{i,j=1}^d h(q) \, g_2^{ij}(q)
   \Big( \frac{\partial}{\partial y_i} \sum_{k=1}^d \lambda_k y^k \Big)(q)
   \Big( \frac{\partial}{\partial y_j} \sum_{l=1}^d \lambda_l y^l \Big)(q)  \\
& = & \sum_{i,j=1}^d h(q) \, g_2^{ij}(q) \, \lambda_i \, \lambda_j
. 
\end{eqnarray*}
But then (\ref{epdrum202;2}) gives
\[
\sum_{i,j=1}^d g_1^{ij}(\tau(q)) \, \lambda_i \, \lambda_j
= \sum_{i,j=1}^d g_2^{ij}(q) \, \lambda_i \, \lambda_j
\]
for all $\lambda_1,\ldots,\lambda_d \in \Ri$ and 
$(g_1^{ij} \circ \tau)(q) = g_2^{ij}(q)$ for all 
$i,j \in \{ 1,\ldots,d \} $.
Hence $g_{1 \, ij}|_{\tau(q)} = g_{2 \, ij}|_q$.
In particular,
\[
g_1|_{\tau(q)}(\frac{\partial}{\partial x^i}, \frac{\partial}{\partial x^j})
= g_2|_q(\frac{\partial}{\partial y^i}, \frac{\partial}{\partial y^j})
= g_2|_q(\tau_* \frac{\partial}{\partial x^i}, \tau_* \frac{\partial}{\partial x^j})
\]
for all $i,j \in \{ 1,\ldots,d \} $.
Hence $\tau$ is an isomorphism from $(M_2,g_2)$ onto $(M_1,g_1)$.\hfill$\Box$

\ruimte

Now the implications \ref{tcdrum101-3}$\Rightarrow$\ref{tcdrum101-1}
and \ref{tcdrum101-4}$\Rightarrow$\ref{tcdrum101-1} in
Theorem~\ref{tcdrum101} follow easily from Lemma~\ref{lcdrum201} and 
Proposition~\ref{pdrum202}.
Substituting $\varphi = \one$ in (\ref{epdrum202;2}) gives 
$\Delta_2 h = 0$ in the proof of Proposition~\ref{pdrum202}.
If $M_2$ is connected this implies that $h$ is constant.
Then the last part in Theorem~\ref{tcdrum101} is obvious.

\subsection*{Acknowledgements}
The second named author is most grateful for the hospitality extended
to him during a fruitful stay at the University of Ulm.
He wishes to thank the University of Ulm for financial support.
Part of this work is supported by the Marsden Fund Council from Government funding, 
administered by the Royal Society of New Zealand.

\end{document}